\documentclass[11pt]{amsart}
\usepackage{amssymb, amsmath, psfrag, epsfig, graphs} 

\newtheorem{thm}{Theorem}[section]
\newtheorem{cor}[thm]{Corollary}
\newtheorem{lem}[thm]{Lemma}
\newtheorem{prop}[thm]{Proposition}
\newtheorem{defn}[thm]{Definition}

\theoremstyle{remark}

\numberwithin{equation}{section}

\newcommand{\de}{\delta}

\newcommand{\Si}{\Sigma}

\newcommand{\Th}{\Theta}

\newcommand{\bD}{\mathbf D}
\newcommand{\CC}{\mathcal C}
\newcommand{\LL}{\mathcal L}

\newcommand{\FF}{\mathcal F}
\renewcommand{\SS}{\mathcal S}

\newcommand{\RR}{\mathcal R}

\newcommand{\Z}{\mathbb Z}
\newcommand{\Q}{\mathbb Q}

\newcommand{\del}{\partial}

\DeclareMathOperator{\Aut}{Aut}
\begin{document}

\title{Sums of lens spaces bounding rational balls} 

\author{Paolo Lisca}
\address{Dipartimento di Matematica ``L. Tonelli''\\ 
Largo Bruno Pontecorvo, 5\\
Universit\`a di Pisa \\
I-56127 Pisa, ITALY} 
\email{lisca@dm.unipi.it}

\date{}

\begin{abstract}
We classify connected sums of three--dimensional lens spaces which smoothly bound rational 
homology balls. We use this result to determine the order of each lens space in the group of 
rational homology 3--spheres up to rational homology cobordisms, and to determine the 
concordance order of each 2--bridge knot. 
\end{abstract}

\maketitle

\section{Introduction}
\label{s:intro}

This paper is a continuation of~\cite{Li}. Both papers can be viewed as providing partial  
answers to the question of Andrew Casson~\cite[Problem~4.5]{Ki} asking which rational homology 3--spheres  
bound rationally acyclic 4--manifolds.  In ~\cite{Li} we identified the set $\RR\subset\Q$ of rational numbers 
$p/q>1$ such that the lens space $L(p,q)$ smoothly bounds a rational homology ball (see Section~\ref{s:coro} for a description of $\RR$). Here we generalize that result by determining which connected sums of lens spaces bound rational balls. Let  
\[
\FF_n := \{\frac{m^2 n}{mnk+1}\ |\ m>k>0,\ (m,k)=1\}\subset\Q,\quad n\geq 2.
\]
\begin{thm}\label{t:main}
Let $Y$ be an oriented 3--manifold homeomorphic to a connected sum of three--dimensional lens 
spaces. Then, $Y$ smoothly bounds a rational homology ball if and only if $Y$ is 
orientation--preserving homeomorphic to a connected sum in which each summand is (possibly 
orientation--reversing) homeomorphic to one of the following oriented 3--manifolds:  
\begin{enumerate}
\item
$L(p,q)$,\ $p/q\in\RR$;
\item
$L(p,q)\# L(p,p-q)$, $p/q>1$;
\item
$L(p_1,q_1)\# L(p_2,q_2)$, $p_i/q_i\in\FF_2$, $i=1,2$;
\item
$L(p,q)\# L(n,n-1)$, $p/q\in\FF_n$ for some $n\geq 2$;
\item
$L(p_1,q_1)\# L(p_2,p_2-q_2)$, $p_i/q_i\in\FF_n$, $i=1,2$, for some $n\geq 2$.
\end{enumerate}
\end{thm}
Let $\Th^3_\Q$ be the group of rational homology 3--spheres up to rational homology cobordism. 
Using Theorem~\ref{t:main} we can determine the order of each lens space viewed as an element 
of the group  $\Th^3_\Q$ of rational homology 3--spheres up to rational homology cobordism. 
Given coprime integers $p>q>0$, denote by $q'$ the unique integer satisfying $p>q'>0$ and 
$qq'\equiv 1\bmod p$. Define 
\[
\SS:=\{p/q\ |\ p>q>0,\ (p,q)=1\ \text{and}\ p=q+q'\}\subset\Q.
\]

\begin{cor}\label{c:lensorder}
The order of the lens space $L(p,q)$ in $\Th^3_\Q$ is:
\begin{itemize}
\item
$1$ if and only if $p/q\in\RR$,
\item
$2$ if and only if $p/q\in(\SS\setminus\RR)\cup\FF_2$,
\item
$\infty$ if and only if $p/q\not\in\RR\cup\SS\cup\FF_2$.
\end{itemize}
\end{cor}
Since the sets $\RR$, $\SS$ and $\FF_2$ have elementary, explicit definitions, 
Corollary~\ref{c:lensorder} reduces the determination of the order of each lens space in $\Th^3_\Q$ 
to an elementary calculation. 

Consistently with~\cite{Li}, in the present paper we define the lens space $L(p,q)$ as the 2--fold cover of $S^3$ branched along the 2--bridge link $K(p,q)$ (see Section~\ref{s:coro} for the definition of $K(p,q)$). 
As a consequence of this fact, Theorem~\ref{t:main} gives information on the relation of link concordance amongst 2--bridge links $K(p,q)\subset S^3$ . Here we content ourselves with deriving a corollary of Theorem~\ref{t:main} in the knot case (i.e.~when $p$ is odd), leaving the case of links to the interested reader. Recall that the~\emph{concordance order} of a knot in $S^3$ is the order of its class 
in the smooth knot concordance group $\CC_1$. The survey paper~\cite{Lv} describes 
what was known about the group $\CC_1$ until 2004. More recent papers contain   
results on knot concordance obtained via Ozsv\'ath--Szab\'o's Heegaard Floer homology~\cite{GRS, JN, OM}. 
In spite of all the efforts made so far, the structure of the group $\CC_1$ remains quite mysterious. For instance, the basic question asking whether $\CC_1$ contains nontrivial elements of finite order different from two is still wide open.  Restricting to 2--bridge knots, results from~\cite{GRS, JN} determine the corresponding concordance orders in finitely many cases. 
\begin{cor}\label{c:2-bridge}
Let $p>q>0$ with $p$ odd. Then, the 2--bridge knot $K(p,q)$ has concordance order:
\begin{itemize}
\item
$1$ if and only if $p/q\in\RR$,
\item
$2$ if and only if $p/q\in\SS\setminus\RR$,
\item
$\infty$ if and only if $p/q\not\in\RR\cup\SS$.
\end{itemize}
Moreover, $K(p,q)$ is simultaneously smoothly slice and amphicheiral if and only if 
\[
(p,q) = \left((2k^2+2k+1)^2,2(2k+1)(k^2+k+1)\right)
\]
for some $k\geq 1$. 
\end{cor}

Corollary~\ref{c:2-bridge} corroborates~\cite[Conjecture~9.4]{Li}, stating that each knot of order two is concordant 
to a negative amphicheiral knot. In fact, since $K(p,q)$ is isotopic to $K(p,q')$ while $K(p,p-q)$ is isotopic to the 
mirror image of $K(p,q)$, by Corollary~\ref{c:2-bridge} each 2--bridge knot of order two is 
amphicheiral~\footnote{Since 2--bridge knots are reversible~\cite{Si}, for such knots the notions of negative and 
positive amphicheirality coincide.}. 

The paper is organized as follows. In Section~\ref{s:coro} we establish 
Corollaries~\ref{c:lensorder} and~\ref{c:2-bridge} assuming Theorem~\ref{t:main}. In Section~\ref{s:sufficiency} we prove that the manifolds (1)--(5) of Theorem~\ref{t:main} bound rational homology balls.  In Section~\ref{s:prelim} we recall some results from~\cite{Li} and prove a few similar results.  In Section~\ref{s:main} we use the results of Section~\ref{s:prelim} to finish the proof of Theorem~\ref{t:main}. 

\section{Proof of Corollaries~\ref{c:lensorder} and~\ref{c:2-bridge}}\label{s:coro}

In this section we prove Corollaries~\ref{c:lensorder} and~\ref{c:2-bridge} assuming Theorem~\ref{t:main}.

\begin{proof}[Proof of Corollary~\ref{c:lensorder}]
By the definition of $\RR$, the lens space $L(p,q)$ has order $1$ in $\Th^3_\Q$, i.e.~it represents 
the trivial class, if and only if $p/q\in\RR$. 

If $L(p,q)\# L(p,q)$ bounds a rational ball, a quick inspection of Theorem~\ref{t:main} shows that 
one of the following must hold:
\[
\text{(1)}\ p/q\in\RR,\quad\text{(2)}\ p=q+q',\quad\text{or (3)}\ p/q\in\FF_2.
\]
On the other hand, if $L(p,q)$ has order $2$ in $\Th^3_\Q$ then (1) does not hold. Hence, by the 
definition of 
$\SS$ we have $p/q\in(\SS\setminus\RR)\cup\FF_2$. 

Conversely, suppose $p/q\in(\SS\setminus\RR)\cup\FF_2$. If $p/q\in\SS$ then $q'=p-q$, 
therefore 
\[
L(p,q)\# L(p,q) = L(p,q)\# L(p,q') =  L(p,q\# L(p,p-q) = L(p,q)\# -L(p,q)
\]
bounds a rational ball, while $L(p,q)$ does not bound such a ball because $p/q\not\in\RR$. 
Therefore $L(p,q)$ has order $2$ in $\Th^3_\Q$. If $p/q\in\FF_2$, by Theorem~\ref{t:main}(3)  
the connected sum $L(p,q)\# L(p,q)$ bounds a rational ball. 

To finish the proof it suffices to show that if $p/q\not\in\RR\cup\SS\cup\FF_2$ then $L(p,q)$ has 
infinite order in $\Th^3_\Q$. But Theorem~\ref{t:main} immediately implies that if $L(p,q)$ has 
finite order then it has either order $1$ or $2$, and 
therefore by what we have already proved $p/q\in\RR\cup\SS\cup\FF_2$.
\end{proof}

\begin{figure}[ht]
\begin{center}
\psfrag{-c1}{${\scriptstyle -c_1}$} 
\psfrag{c2}{${\scriptstyle c_2}$}
\psfrag{-c3}{${\scriptstyle -c_3}$} 
\psfrag{-cn-2}{${\scriptstyle -c_{n-2}}$} 
\psfrag{cn-1}{${\scriptstyle c_{n-1}}$}
\psfrag{-cn-1}{${\scriptstyle -c_{n-1}}$}
\psfrag{cn}{${\scriptstyle c_n}$} 
\psfrag{-cn}{${\scriptstyle -c_n}$} 
\psfrag{n even}{($n$ even)}
\psfrag{n odd}{($n$ odd)}
\psfrag{s}{\small $s$}
\psfrag{crossings}{$|s|$ crossings} 
\psfrag{s>0}{if $s\geq 0$}
\psfrag{s<0}{if $s\leq 0$} 
\includegraphics[height=8cm]{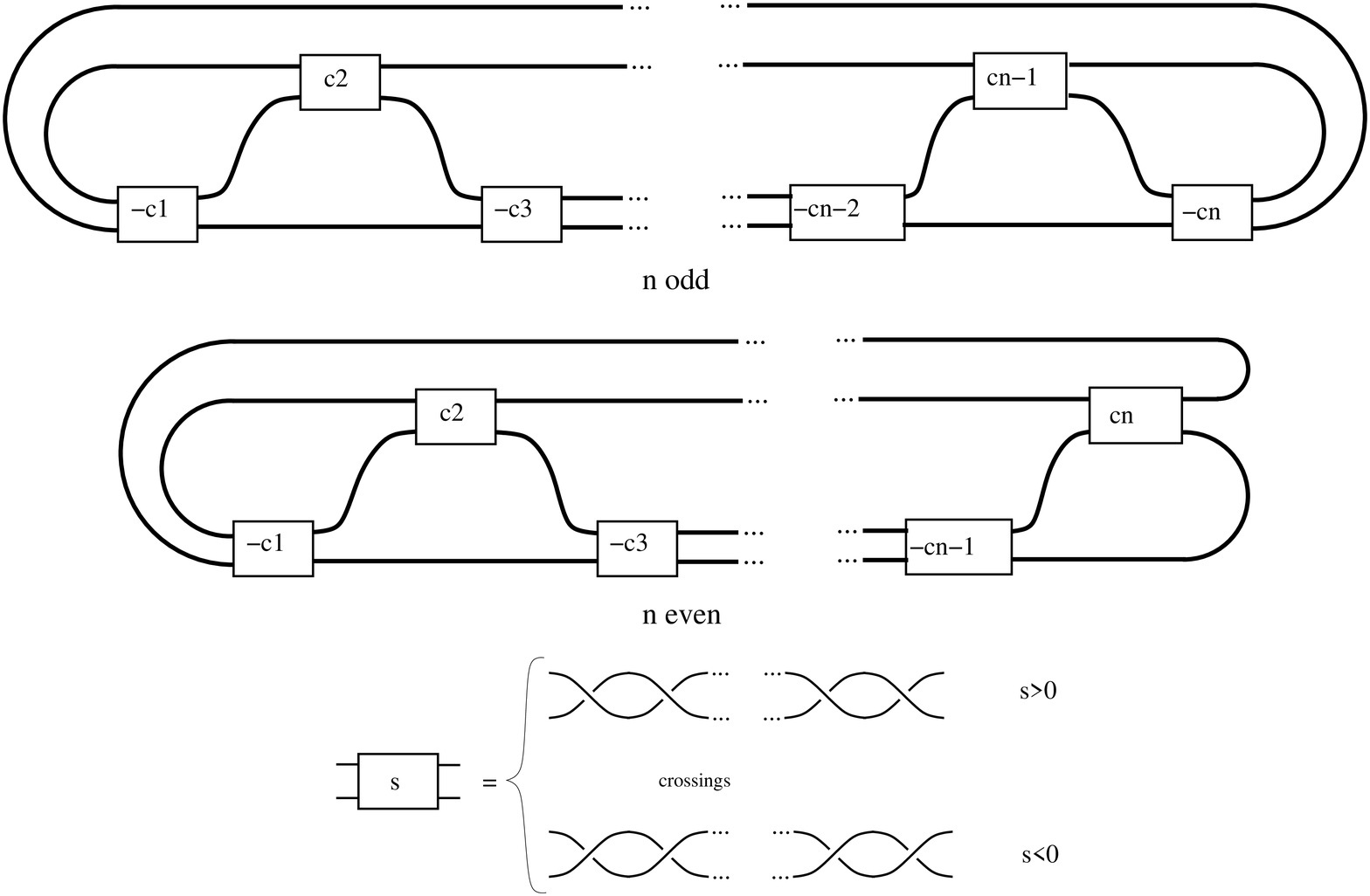}
\end{center}
\caption{The $2$--bridge link $\LL(c_1,\ldots,c_n)$}
\label{f:fig1}
\end{figure}

Before proving  Corollary~\ref{c:2-bridge} we need to recall the definitions of the links $K(p,q)$ and the set $\RR$.
 A link in $S^3$ is called~\emph{$2$--bridge} if it can be isotoped
so it has exactly two local maxima with respect to a standard
height function. Figure~\ref{f:fig1} represents the $2$--bridge link
$\LL(c_1,\ldots,c_n)$, where $c_i\in\Z$, $i=1,\ldots,n$.
Given coprime integers $p>q>0$, we can alway write $p/q=[c_1,\ldots,c_n]^+$, 
$c_i>0$, $i=1,\ldots, n$. The 2--bridge link $K(p,q)$ is, by definition, $\LL(c_1,\ldots,c_n)$. When 
$p$ is even, $K(p,q)$ is a 2--component link, when $p$ is odd $K(p,q)$ is a knot. It is 
well--known~\cite[Chapter~12]{BZ} that $K(p,q)$ and $K(\overline p,\overline q)$ are isotopic if and only if $p=\overline p$ and either $q = \overline q$ or $q{\overline q}\equiv 1\pmod p$, and that every $2$--bridge link is isotopic to some $K(p,q)$. Moreover, $K(p,p-q)$ is isotopic to the mirror image of $K(p,q)$. 

Recall~\cite{Li} that a rational number $p/q>1$ with $(p,q)=1$ belongs to $\RR$ if and only 
if (i) $p=m^2$ and (ii) $q$, $p-q$ or $q'$ is of one of the following types:
\begin{enumerate}
\item[(a)]
$mk\pm 1$ with $m>k>0$ and $(m,k)=1$;
\item[(b)]
$d(m\pm 1)$, where $d>1$ divides $2m\mp 1$;
\item[(c)]
$d(m\pm 1)$, where $d>1$ is odd and divides $m\pm 1$.
\end{enumerate}

\begin{proof}[Proof of Corollary~\ref{c:2-bridge}]
If $K(p,q)$ has finite order in $\CC_1$ then $L(p,q)$ has finite order in $\Th^3_\Q$. Therefore by 
Corollary~\ref{c:lensorder} 
we have $p/q\in\RR\cup\SS\cup\FF_2$. If $p/q\in\FF_2$ then $p$ would be even, while in our 
case $p$ is odd because we are 
assuming that $K(p,q)$ is a knot. Hence, $p/q\in\RR\cup\SS$. Moreover, $K(p,q)$ is smoothly slice 
if and only if $p/q\in\RR$ by~\cite[Corollary~1.3]{Li}, while if $p/q\in\SS$ then $q'=p-q$ and 
therefore $K(p,q)=K(p,q')=K(p,p-q)$, so 
\[
K(p,q)\# K(p,q) = K(p,q\# K(p,p-q) = K(p,q)\# -K(p,q)
\]
is smoothly slice. This shows that $K(p,q)$ has finite order in $\CC_1$ if and only if either (i) 
$p/q\in\RR$ and $K(p,q)$ has order $1$ or (ii) $p/q\in\SS\setminus\RR$ and $K(p,q)$ has order 
$2$. This concludes the proof of the part of the statement about concordance orders. 

For the last statement, in Case (a) of the definition of $\RR$ it is easy to check that $q'=m(m-k)\pm 1$, hence $q+q'=m^2\pm 2\neq p$.  In Case (b), if $q=d(m-1)$ then $q'=h(m-1)$, where $h$ is defined by $dh=2m+1$. Then, $q+q'=(d+h)(m-1)=m^2=p$ is impossible. If $q=d(m+1)$ we have $q'=m^2-h(m+1)$, with $dh=2m-1$. Thus $q+q'=p$ if and only if $d=h$, which implies $m=(d^2+1)/2$, so $(p,q)$ is of the form $((d^2+1)^2/4,d(d^2+3)/2)$. Observing that $d>1$ is odd and setting $d=2k+1$ one obtains the statement of the lemma. In Case (c), if $q=d(m-1)$ then $q'=h(2m+1)$, with  $dh=m+1$, and the equality $q+q'=p$ gives a contradiction as in Case (b). If $q=d(m+1)$ then $q'=m^2-h(2m-1)$, and the condition $q+q'=p$ implies $m=(d^2+1)/2$ as before, to the effect that $(p,q)$ must be of the stated form. 
\end{proof}

\section{Proof of Theorem~\ref{t:main}: first half}\label{s:sufficiency}

The purpose of this section is to prove the first half of Theorem~\ref{t:main}, i.e.~the following
\begin{prop}\label{p:firsthalf}
Each of the manifolds (1)--(5) listed in Theorem~\ref{t:main} bounds a rational homology ball.
\end{prop}

We start with some preparation. Given integers $a_1,\ldots,a_h\geq 1$ and $b_1,\ldots,b_k\geq 2$, let 
\[
[a_1,\ldots, a_k]^+ := a_1 + \cfrac{1}{a_2 +
       \cfrac{1}{\ddots +
        \cfrac{1}{a_h}
}},\ 
[b_1,\ldots, b_k]^- := b_1 - \cfrac{1}{b_2 -
       \cfrac{1}{\ddots -
        \cfrac{1}{b_k}
}}.
\]
For any integer $t\geq 0$ we shall write 
\begin{equation*}
(\ldots,2^{[t]},\ldots) := (\ldots,\overbrace{2,\ldots,2}^t,\ldots).
\end{equation*}
For example,
\[
(3,2^{[3]},4) = (3,2,2,2,4).
\]
Let $a_1,\ldots,a_{2d}$ be positive integers and $d\geq 2$. Then, the following 
identity holds~\cite[Proposition~2.3]{PP}:
\begin{equation}\label{e:cfrac}
[a_1,\ldots,a_{2d}]^+ = 
[a_1+1,2^{[a_2-1]},a_3+2,2^{[a_4-1]},\ldots,a_{2d-1}+2,2^{[a_{2d}-1]}]^-
\end{equation}
Also, recall (see e.g.~\cite[Appendix]{OW}) that if $p/q>1$ then 
\[
\frac pq=[a_1,a_2,\ldots,a_{h-1},a_h]^- 
\]
for some $a_i\geq 2$ and  
\begin{equation}\label{e:fracreverse}
\frac p{q'} = [a_h,a_{h-1},\ldots,a_2,a_1]^- ,
\end{equation}
where $p>q'>0$ and $qq'\equiv 1\pmod p$. 

We define the {\bf reverse} of a string $(a_1,\ldots,a_k)$ to be  
$(a_k,\ldots,a_1)$, and the  {\bf negative} (respectively {\bf positive}) 
{\bf fraction} of  $(a_1,\ldots,a_k)$ to be the rational number 
$[a_1,\ldots, a_k]^-$ (respectively $ [a_1,\ldots,a_k]^+$). 

\begin{lem}\label{l:fraction}
Suppose that a string $\SS$ is obtained from $(2,n+1,2)$, $n\geq 2$, by a finite 
sequence of operations each of which is of one of the following types: 
\begin{enumerate}
\item
$(n_1,n_2,\ldots,n_{b-1},n_b)\mapsto (2,n_1,n_2,\ldots,n_{b-1},n_b+1)$,
\item
$(n_1,n_2,\ldots,n_{b-1},n_b)\mapsto (n_1+1,n_2,\ldots,n_{b-1},n_b,2)$.
\end{enumerate}
Then, the negative fraction of $\SS$ is of the form  
\[
\frac{m^2 n}{mnk+1},\quad m>k>0,\quad (m,k)=1.
\]
Moreover, the negative fraction of either $\SS$ or the reverse of $\SS$ is equal 
to 
\begin{equation}\label{e:posexp}
[c_s,c_{s-1},\ldots,c_2,c_1,1,n-1,c_1+1,c_2,\ldots,c_s]^+
\end{equation}
for some integers $c_1,\ldots,c_s\geq 1$, where for $s=1$ the above 
formula is to be interpreted as $[c_1,1,n-1,c_1+1]^+$.
\end{lem}

\begin{proof}
The negative fraction of $(2,n+1,2)$ is $4n/(2n+1)$, which is 
of the form $m^2 n/(mnk+1)$ for $m=2$ and $k=1$. Therefore it 
suffices to show that if the negative fraction of a string $\SS$ is of the form 
$m^2 n/(mnk+1)$, then the negative fractions of the strings obtained from $\SS$ by 
applying the operations (1) and (2) are of the same form. Let 
\begin{align}\label{e:qprimes}
p:=m^2 n,\quad q:=mnk+1,\quad q':=mn(m-k)+1,\\ 
t:=nk(m-k)+1,\quad q'':=2q-nk^2=q+t.
\end{align}
Then, we have $0<q'<p$, $(p,q')=1$,
\[
q q' = 1 + tp \quad\text{and}\quad q'' q' = 1 + t(p+q').
\]
Therefore if $p/q=[a_1,\ldots,a_h]^-$, by Equation~\eqref{e:fracreverse}
\[
[a_h+1,a_{h-1},\ldots,a_2,a_1]^- = \frac{p+q'}{q'}.
\]
Since $0<q''=q+t<q+q'<p+q'$ and $(q'',p+q')=1$, we have 
\[
[a_1,a_2,\ldots,a_{h-1},a_h+1]^- = \frac{p+q'}{q''}
\]
and therefore 
\begin{equation}\label{e:F}
[2,a_1,a_2,\ldots,a_{h-1},a_h+1]^- = 
2 - \frac{q''}{p+q'} = \frac{(2m-k)^2 n}{(2m-k)nm + 1},
\end{equation}
with $2m-k>m$ and $(2m-k,m)=1$. Applying Equation~\eqref{e:F} to
\[
[a_h,\ldots,a_1]^- = \frac{m^2 n}{mn(m-k)+1}
\]
gives 
\[
[2,a_h,\ldots,a_1+1]^- = \frac{(m+k)^2 n}{(m+k)nm+1}.
\]
By Equation~\eqref{e:fracreverse} and the formula for $q'$ in Equation~\eqref{e:qprimes} 
we have 
\begin{equation}\label{e:induction}
[a_1+1,a_2,\ldots,a_{h-1},a_h,2]^- = 
\frac{(m+k)^2 n}{(m+k)nk + 1}.
\end{equation}
This proves the first part of the lemma. To prove the second part, observe that 
a straightforward induction using the definition shows that either $\SS$ or 
its reverse is of the form 
\[
\begin{cases}
{\scriptstyle (c_s+1,2^{[c_{s-1}-1]},c_{s-2}+2,\ldots,c_3+2,2^{[c_2-1]},c_1+2,
n+1, 2^{[c_1]},c_2+2,\ldots,c_{s-1}+2,2^{[c_s-1]})}\quad\text{or}\\ 
(c_1+1,n+1,2^{[c_1]}).
\end{cases}
\]
for some positive integers $c_1,\ldots,c_s$. Applying Formula~\eqref{e:cfrac} we 
obtain Expansion~\eqref{e:posexp}.
\end{proof}

\begin{lem}\label{l:expansion}
Let $m>k>0$ and $n\geq 1$ be integers, with $(m,k)=1$. Then, 
\[
\frac{m^2n}{mnk+1} = [a_1,\ldots,a_h]^-,
\]
where the string $(a_1,\ldots,a_h)$ is obtained from $(2,n+1,2)$ by a sequence 
of operations of types (1) or (2) as in Lemma~\ref{l:fraction}. 
\end{lem}

\begin{proof}
We argue by induction on $m\geq 2$. For $m=2$ we have $k=1$ and 
\[
\frac{4n}{2n+1} = [2,n+1,2]^-.
\]
Now assume $m>2$ and the conclusion of the lemma for any fraction 
\[
\frac{p^2n}{pnl+1},\quad
p>l>0,\quad (p,l)=1,
\]
with $p<m$. If $m=2k$ then $m=2$, therefore 
we have either $m>2k$ or $m<2k$. Since 
\[
(mnk+1)(mn(m-k)+1)\equiv 1\pmod{m^2n},
\]
by Equation~\eqref{e:fracreverse}
\[
\frac{m^2n}{mnk+1} = [a_1,\ldots,a_h]^- \ \Longleftrightarrow\ 
\frac{m^2n}{mn(m-k)+1} = [a_h,\ldots,a_1]^-.
\]
Since string reversal turns the operations (1) and (2) of Lemma~\ref{l:fraction} into each other, 
up to replacing $k$ with $m-k$ we may assume $m>2k$ without losing generality. 
Setting $p:=m-k$ and $l:=k$ we have $m>p$, $p>l>0$ and $(p,l)=1$. 
By the induction hypothesis 
\[
\frac{p^2n}{pnl+1} = [b_1,\ldots,b_q]^-, 
\]
where $(b_1,\ldots,b_q)$ is obtained from $(2,n+1,2)$ by a sequence 
of operations of types (1) or (2) as in Lemma~\ref{l:fraction}. Then, by 
Equation~\eqref{e:induction} we have 
\[
(a_1,\ldots,a_h)=(b_1+1,\ldots,b_q,2).
\]
This concludes the proof. 
\end{proof}

\begin{lem}\label{l:move}
Let $p/q\in\FF_n$, $n\geq 2$. Then, there exists a ribbon move which turns the link $K(p,q)$ 
into a split link given by the union of $K(n,1)$ and an unknot. 
\end{lem}

\begin{figure}[ht]
\begin{center}
\psfrag{-c1}{$-c_1$} 
\psfrag{c1+1}{$c_1+1$} 
\psfrag{c2}{${c_2}$}
\psfrag{-c2}{$-c_2$}
\psfrag{1-n}{${\scriptstyle 1-n}$} 
\psfrag{cs}{$c_s$}
\psfrag{-cs}{$-c_s$}
\psfrag{-n}{$-n$} 
\psfrag{iso}{isotopy}
\includegraphics[height=15cm]{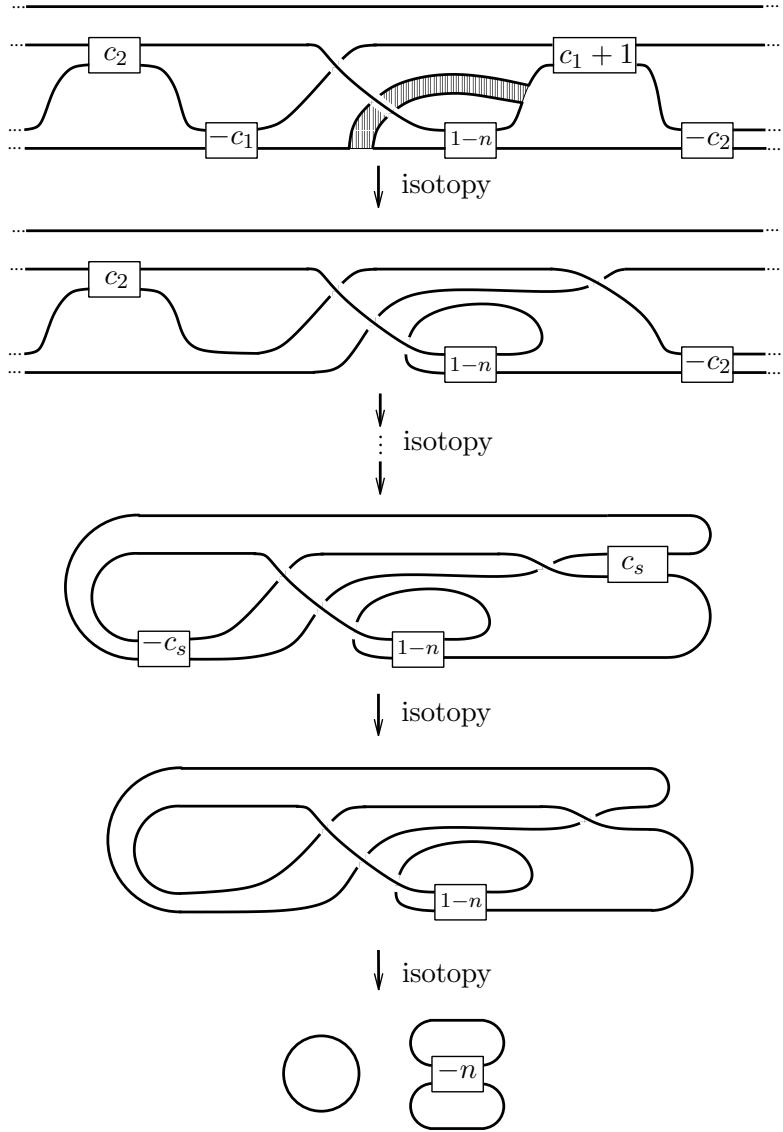}
\end{center}
\caption{Ribbon move on $K(m^2 n,mnk+1)$ and isotopies}
\label{f:fig3}
\end{figure}

\begin{proof}
By Lemmas~\ref{l:fraction} and~\ref{l:expansion}, up to isotopy we have the equality
\[
K(m^2n,mnk+1) = \LL(c_s,\ldots,c_1,1,n-1,c_1+1,c_2,\ldots,c_s)
\]
for some integers $c_1,\ldots,c_s\geq 1$. The rest of the proof consists of  
Figure~\ref{f:fig3}.
\end{proof}	

The following lemma can be extracted from~\cite[Proof of Theorem~1.2]{Li}. We include the proof  
for the reader's convenience. 

\begin{lem}\label{l:bounding}
Let $Y$ be an oriented 3--manifold given as the 2--fold cover of $S^3$ branched along a link $L$ 
which bounds a ribbon surface $\Si\looparrowright S^3$ with $\chi(\Si)=1$. Then, $Y$ smoothly bounds 
a rational homology 4--ball. 
\end{lem}

\begin{proof}
The 4--manifold $W$ given as the 2--fold cover of the 4--ball $B^4$ branched along a pushed--in copy 
of $\Si$ smoothly bounds $Y$. Moreover, $W$ has a handle decomposition with only 0--, 1-- and 
2--handles (see e.g.~\cite[lemma at pages 30--31]{CH}). Therefore, from
\[
b_0(W)-b_1(W)+b_2(W) = \chi(W) = 2\chi(B^4) - \chi(\widetilde\Si) = 1
\]
we deduce $b_1(W)=b_2(W)$. On the other hand, since $b_1(\del W)=0$
and $H_1(W,\del W;\Q)\cong H^3(W;\Q)=0$, the homology exact sequence
of the pair $(W,\del W)$ gives $b_1(W)=0$, so it follows that
$H_*(W;\Q)\cong H_*(B^4;\Q)$
\end{proof}

\begin{proof}[Proof of Proposition~\ref{p:firsthalf}]
In Case (1) the proposition follows from the results of~\cite{Si} (see~\cite{Li} for complete 
details). In Case (2) the result is well--known. To treat Cases (3)--(5), we will apply 
Lemma~\ref{l:bounding}. First observe that, as shown in Figure~\ref{f:fig4}, there is a connected 
sum $K(n,1)\# K(n,n-1)$ which can be turned into an unlink of 2 components by a ribbon move. 
\begin{figure}[ht]
\begin{center}
\psfrag{-n}{$-n$} 
\psfrag{n}{$n$}
\psfrag{=}{$=$}
\psfrag{iso}{isotopy}
\includegraphics[width=10cm]{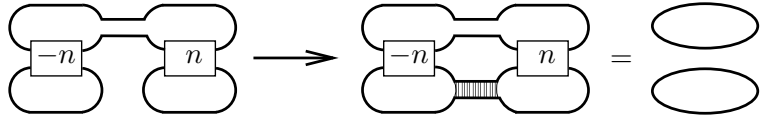}
\end{center}
\caption{Ribbon move on $K(n,1)\# K(n,n-1)$}
\label{f:fig4}
\end{figure}
 On the other hand, applying Lemma~\ref{l:move} we see that if ${p_i}/{q_i}\in\FF_2$ then there is a connected sum $K(p_1,q_1)\# K(p_2,q_2)$ which can be turned into the connected sum $K(2,1)\# K(2,1)$ of Figure~\ref{f:fig4} (for $n=2$) by 2 ribbon moves. 
Therefore, a connected sum $K(p_1,q_1)\# K(p_2,q_2)$ with ${p_i}/{q_i}\in\FF_2$ can be turned into a 4--component unlink by 3 ribbon moves. So we conclude that such a link bounds a ribbon surface $\Si$ with $\chi(\Si)=1$. Exactly the same argument, with the same conclusion, applies to a connected sum $K(p_1,q_1)\# K(p_2,p_2-q_2)$ with ${p_i}/{q_i}\in\FF_n$, $n\geq 2$. Since the 2--fold cover of $S^3$ branched along any connected sum $K(p_1,q_1)\# K(p_2,q_2)$ is diffeomorphic to $L(p_1,q_1)\# L(p_2,q_2)$, by  Lemma~\ref{l:bounding} this suffices to prove the proposition in Cases (3) and (5). A similar construction shows that a connected sum $K(p,q)\# K(n,n-1)$ can be reduced to a 3--component unlink by 2 ribbon moves. Therefore the same  connected sum also bounds a ribbon surface $\Si$ with $\chi(\Si)=1$, and applying Lemma~\ref{l:bounding} as before this proves the proposition in Case (4). 
\end{proof}

\section{Algebraic interlude}\label{s:prelim}

In this section we recall some definitions and results from~\cite{Li} and 
we establish some new results having a similar flavor. This material will be used in 
Section~\ref{s:main}. 

\part*{Recollection of previous results}

Let $\bD$ denote the intersection lattice $(\Z,(-1))$, and 
$\bD^n$ the orthogonal direct sum of $n$ copies of $\bD$. Elements of $\bD^n$ will also 
be called~\emph{vectors}. Fix generators $e_1,\ldots,e_n\in\bD^n$ such that
\[
e_i\cdot e_j=-\de_{ij},\quad i,j=1,\ldots,n. 
\]

Given a subset
$S=\{v_1,\ldots,v_n\}\subseteq\bD^n$, define
\[
E^S_i:=\{j\in\{1,\ldots,n\}\ |\ v_j\cdot e_i\neq 0\},\quad
i=1,\ldots,n,
\]
\[
V_i:=\{j\in\{1,\ldots, n\}\ |\ e_j\cdot v_i\neq 0\},\quad
i=1,\ldots,n,
\]
and 
\[
p_i(S):=|\{j\in\{1,\ldots,n\}\ |\ |E^S_j|=i\}|,\quad
i=1,\ldots,n.
\]

Let $v_1,\ldots, v_n\in\bD^n$ be elements such that, for
$i,j\in\{1,\ldots,n\}$,
\begin{equation}\label{e:conds}
v_i\cdot v_j = 
\begin{cases}
-a_i\leq -2\quad& \text{if}\ i=j,\\
0\ \text{or}\ 1 \quad& \text{if}\ |i-j|=1,\\
0\quad& \text{if}\ |i-j|>1,
\end{cases}
\end{equation}
for some integers $a_i$, $i=1,\ldots,n$.

Let $S=\{v_1,\ldots, v_n\}\subseteq\bD^n$ be a subset which 
satisfies~\eqref{e:conds}. Define the {\bf intersection graph} of $S$ as the graph 
having as vertices $v_1,\ldots,v_n$, and an edge between $v_i$ and $v_j$ if and only 
if $v_i\cdot v_j=1$ for $i, j=1,\ldots,n$. The number of connected components of 
the intersection graph of $S$ will be denoted by $c(S)$. An element $v_j\in S$ 
is {\bf isolated}, {\bf final} or {\bf internal} if the quantity
\[
\sum^n _{\substack{i=1\\ i\neq j}} 
(v_i\cdot v_j) 
\]
is equal to, respectively, $0$, $1$ or $2$. In other words, $v_j$ is isolated or 
final if it is, respectively, an isolated vertex or a leaf of the intersection graph, 
and it is internal otherwise. Two elements $v, w\in\bD^n$ are {\bf linked} if there 
exists $e\in\bD^n$ with $e\cdot e=-1$ such that 
\[
v\cdot e\neq 0,\quad\text{and}\quad w\cdot e\neq 0.
\]
A set $S\subseteq\bD^n$ is~{\bf irreducible} if, given two elements
$v,w\in S$, there exists a finite sequence of elements of $S$
\[
v_0=v, v_1,\ldots, v_k=w, 
\]
such that $v_i$ and $v_{i+1}$ are linked for  $i=0,\ldots, k-1$. 
A set which is not irreducible is~{\bf reducible}.

A subset $S=\{v_1,\ldots,v_n\}\subseteq\bD^n$ is~{\bf good} if it is
irreducible and its elements satisfy~\eqref{e:conds}. If moreover $v_i\cdot v_j=1$ 
whenever $|i-j|=1$, we say that $S$ is~{\bf standard}. 

Given a subset $S=\{v_1\ldots,v_n\}\subseteq\bD^n$, define
\[
I(S):=\sum_{i=1}^n (-v_i\cdot v_i -3) \in \Z.
\]
Let $\Aut(\bD^n)$ be the group of automorphisms of $\bD^n$ as an intersection lattice.

\begin{lem}[\cite{Li}, Lemma~2.4]\label{l:n=3}
Let $S=\{v_1,v_2,v_3\}\subseteq\bD^3=\langle e_1,e_2,e_3\rangle$ be a
good subset with $I(S)<0$. Then, up to applying to $S$ an element of
$\Aut(\bD^3)$ and possibly replacing $(v_1, v_2,v_3)$ with
$(v_3,v_2,v_1)$, one of the following holds:
\begin{enumerate}
\item
$(v_1,v_2,v_3)=(e_1-e_2,e_2-e_3,-e_2-e_1)$,
\item
$(v_1, v_2,v_3)=(e_1-e_2,e_2-e_3,e_1+e_2+e_3)$,
\item
$(v_1,v_2,v_3)=(e_1+e_2+e_3,-e_1-e_2+e_3,e_1-e_2)$.
\end{enumerate}
Moreover, 
\[
(p_1(S),p_2(S),c(S),I(S))=
\begin{cases}
(1,1,1,-3)\quad\text{in case $(1)$},\\
(0,2,2,-2)\quad\text{in case $(2)$},\\
(0,1,2,-1)\quad\text{in case $(3)$}.
\end{cases}
\]
In particular, $(a_1,a_2,a_3)\in\{(2,2,2),(2,2,3),(3,3,2)\}$.
\qed\end{lem}

\begin{lem}[\cite{Li}, Lemma~2.6]\label{l:negsum}
Let $p>q\geq 1$ be coprime integers, and suppose that 
\[
\frac pq = [a_1,\ldots,a_n]^-,\quad
\frac p{p-q} = [b_1,\ldots,b_m]^-,
\]
with $a_1,\ldots,a_n\geq 2$ and $b_1,\ldots,b_m\geq 2$. Then, 
\[
\sum_{i=1}^n (a_i-3) + \sum_{j=1}^m (b_j-3) = -2.
\]
\qed\end{lem}

Given elements $e, v\in\bD^n$ with $e\cdot e = -1$, we denote by
$\pi_e(v)$ the projection of $v$ in the direction orthogonal to $e$:
\[
\pi_e(v):=v+(v\cdot e) e\in\bD^n.
\]

\begin{defn}\label{d:contraction}
Let $S=\{v_1,\ldots,v_n\}\subseteq\bD^n$ be a subset
satisfying~\eqref{e:conds} and such that $|v_i\cdot e_j|\leq 1$ for
every $i,j=1,\ldots,n$. Suppose that there exist $1\leq h,s,t\leq n$
such that
\[
E_h^S=\{s,t\}\quad\text{and}\quad a_t>2.
\]
Then, we say that the subset $S'\subseteq\langle e_1,\ldots,e_{h-1},
e_{h+1},\ldots, e_n\rangle\cong\bD^{n-1}$ defined by
\[
S':=S\setminus\{v_s,v_t\}\cup\{\pi_{e_h}(v_t)\}
\]
is obtained from $S$ by a {\bf contraction}, and we write $S\searrow
S'$. We also say that $S$ is obtained from $S'$ by
an {\bf expansion}, and we write $S'\nearrow S$.
\end{defn}

\begin{defn}\label{d:bad}
Let $S'=\{v_1,\ldots, v_n\}\subseteq\bD^n$, $n\geq 3$, be a good
subset, and suppose there exists $1<s<n$ such that
$C'=\{v_{s-1},v_s,v_{s+1}\}\subseteq S'$ is a connected component
of the intersection graph of $S'$, with 
$v_{s-1}\cdot v_{s-1}=v_{s+1}\cdot v_{s+1}= -2$, $v_s\cdot v_s < -2$
and $E^{S'}_j=\{s-1,s,s+1\}$ for some $j$. 
Let $S\subseteq\bD^m$ be a subset of cardinality $m\geq n$ obtained
from $S'$ by a sequence of expansions by final $(-2)$--vectors
attached to $C'$, so that $c(S)=c(S')$ and there is a natural 
one--to--one correspondence between the sets of connected components of the
intersection graphs of $S$ and $S'$. Then, the connected component
$C\subseteq S$ corresponding to $C'\subseteq S'$ is a~{\bf bad
component} of $S$.  The number of bad components of $S$ will be
denoted by $b(S)$.
\end{defn}

In the arguments of the next two sections we shall need the following 
Proposition~\ref{p:I<0} and Lemmas~\ref{l:extracted} and~\ref{l:I=-2}. 

\begin{prop}[\cite{Li}, Corollary~5.4]\label{p:I<0}
Suppose that $n\geq 3$, and let $S_n=\{v_1,\ldots,v_n\}\subseteq\bD^n$
be a good subset with no bad components and such that $I(S_n) <
0$. Then $I(S_n)\in\{-1,-2,-3\}$, there exists a sequence of
contractions
\begin{equation}\label{e:expans}
S_n\searrow S_{n-1}\searrow\cdots\searrow S_3
\end{equation}
such that, for each $k=3,\dotsc,n-1$ the set $S_k$ is good, has no
bad components and we have either
\begin{equation}\label{e:Ic1}
(I(S_k),c(S_k))=(I(S_{k+1}),c(S_{k+1}))
\end{equation}
or
\begin{equation}\label{e:Ic2}
I(S_k)\leq I(S_{k+1})-1\quad\text{and}\quad c(S_k)\leq c(S_{k+1})+1.
\end{equation}
Moreover:
\begin{enumerate}
\item
If $p_1(S_n)>0$ then $I(S_n)=-3$, $S_n$ is standard and one can choose
the above sequence so that $I(S_k)=-3$ and $S_k$ is standard for every
$k=3,\ldots, n-1$.
\item
If $I(S_n)+c(S_n)\leq 0$ then $S_3$ is given, up to applying an
automorphism of $\bD^3$, by either (1) or (2) in Lemma~\ref{l:n=3}; if
$I(S_n)+c(S_n)<0$ then the former case occurs.
\end{enumerate}
\qed\end{prop}

\begin{lem}[\cite{Li}, Lemma~6.3]\label{l:extracted}
Let $S_3\subset \bD^3$ be a good subset with $I(S_3)=-3$ and $c(S_3)=1$.
Suppose that $S_3\nearrow\cdots\nearrow S_k$ is a sequence of expansions
such that, for each $h=3,\ldots,k$, $S_h$ is good, has no bad component and
$(I(S_h),c(S_h))=(-3,1)$. Then, it is not possible to expand $S_k$ by 
an isolated $(-3)$--vector. 
\qed\end{lem}

\begin{lem}[\cite{Li}, Lemma~6.2]\label{l:I=-2}
Let $S_3\nearrow\cdots\nearrow S_n$ be a sequence of expansions such
that, for each $k=3,\ldots,n$, $S_k$ is good, has no bad component and
$(I(S_k),c(S_k))=(-2,2)$. Then,
\begin{enumerate}
\item
$p_1(S_n)=0$, $p_2(S_n)=2$ and $p_3(S_n)=n-2$. 
\item
If $E^{S_n}_i=\{t,t'\}$ then $v_t$ and $v_{t'}$ are not internal and 
exactly one of them has square $-2$. 
\item
If $v_t\in S_n$ is not internal then there exists $i\in V_t$ such that 
$|E^{S_n}_i|=2$.
\end{enumerate}
\qed\end{lem}

\part*{New results}

\begin{lem}\label{l:aux}
Let $S_3\nearrow\cdots\nearrow S_n$ be a sequence of expansions of good 
subsets without bad components such that $(I(S_k),c(S_k))=(-2,2)$ for 
every $k=3,\ldots,n$. Then, it is not possible to expand $S_n$ by an isolated 
$(-3)$--vector. 
\end{lem}

\begin{proof} 
Let $S_n=\{{\overline v_1},\ldots, {\overline v_n}\}
\subset\langle e_1,\ldots, e_n\rangle\cong\bD^n$, and suppose that 
\[
\bD^{n+1}\cong\langle e_1,\ldots,e_{n+1}\rangle\supset
S_{n+1} = \{v_1,\ldots,v_{n+1}\}
\searrow S_n
\] 
is a contraction obtained by eliminating the isolated vector $v_{n+1}\in S_{n+1}$ of square $-3$. 
Since $S_n$ is good, has no 
bad components and $I(S_n)<0$, by~\cite[Proposition~5.2]{Li} we have 
$|{\overline v_i}\cdot e_j|\leq 1$ for every $i,j=1,\ldots, n$. Moreover, we may assume without loss 
$v_{n+1}=e_1+e_2+e_{n+1}$, with $|E_{n+1}^{S_{n+1}}|=2$. Then, since 
\begin{multline*}
|E_1^{S_n}| + |E_2^{S_n}| + 1 = 
|E_1^{S_{n+1}\setminus\{v_{n+1}\}}| + |E_2^{S_{n+1}\setminus\{v_{n+1}\}}| + 
|E_{n+1}^{S_{n+1}\setminus\{v_{n+1}\}}| \\ 
\equiv 
\sum_{i\neq n+1} v_i\cdot v_{n+1} \pmod 2 = 0 \pmod 2,
\end{multline*}
we see that $|E_1^{S_n}| + |E_2^{S_n}|$ is odd. Therefore, in view of  
Lemma~\ref{l:I=-2}(1) we may assume $|E_1^{S_n}|=2$ and $|E_2^{S_n}|=3$.
We claim that this is impossible. In fact, observe that by Lemma~\ref{l:I=-2}(2)
we have $E_1^{S_n}=\{t, t'\}$, with ${\overline v_t}$ and ${\overline v_{t'}}$ not internal and 
${\overline v_t}\cdot {\overline v_t}= -2$, say. Moreover, it is easy to check that ${\overline v_t}$ 
is final and $v_t\cdot e_{n+1}=0$. Therefore, since $v_{n+1}\cdot v_t = 0$ we have 
$v_t = \pm (e_1-e_2)$. Since $E_2^{S_n} = \{t,s,r\}$ for some $s,r$ and ${\overline v_t}$ is 
not internal, we can assume ${\overline v_t}\cdot {\overline v_s}=0$. This implies 
$e_1\cdot {\overline v_s}\neq 0$, therefore ${\overline v_s}={\overline v_{t'}}$. 
It follows that $e_1\cdot {\overline v_r}=0$, hence ${\overline v_r}\cdot {\overline v_t}=1$. 
We conclude ${\overline v_t}\cdot {\overline v_{t'}}=0$, so $v_{t'} = \pm (e_1+e_2) + \cdots$, 
which is incompatible with the fact that $v_{n+1} = e_1 + e_2 + e_{n+1}$ is isolated in 
$S_{n+1}$. 
\end{proof}

\begin{lem}\label{l:cleq2(1)}
Let $S\subseteq\bD^n$ be a good subset with no bad components and  
$I(S)<0$. Then, $c(S)\leq 2$. 
\end{lem}

\begin{proof}
If $n\leq 2$ the conclusion is obvious, so we assume $n\geq 3$. 
By contradiction, suppose $c(S)\geq 3$. In view of Proposition~\ref{p:I<0} we have 
$I(S)\in\{-1,-2,-3\}$ and there is a sequence of contractions 
\begin{equation}\label{e:contr}
S=:S_n\searrow S_{n-1}\searrow\cdots\searrow S_3\subseteq\bD^3,
\end{equation}
where for each $k=3,\ldots, n$ the subset $S_k$ is  good, has no bad 
components and either~\eqref{e:Ic1} or~\eqref{e:Ic2} holds. In particular, 
$I(S_3)\leq I(S_4)\leq\cdots\leq I(S_n)$. 

{\em First case: $I(S_n)=-3$.} Since by Lemma~\ref{l:n=3} we have $-3\leq I(S_3)$,  
in this case $I(S_k)=-3$ for every $k$. Hence~\eqref{e:Ic1} holds for every $k$, 
implying $c(S_3)=c(S_n)\geq 3$, and contradicting Lemma~\ref{l:n=3}. 

{\em Second case: $I(S_n)=-2$.} If $I(S_3)=-2$ then arguing as in the previous 
paragraph one gets a contradiction, so we may assume $I(S_3)=-3$. In this 
case~\eqref{e:Ic1} must hold for all but at most one index $k$. This clearly 
implies $c(S_3)\geq c(S_n)-1\geq 2$, which contradicts Lemma~\ref{l:n=3}.

{\em Third case: $I(S_n)=-1$.} If $I(S_3)=-1$ one gets a contradiction as before. If 
$I(S_3)=-2$ then in Sequence~\eqref{e:contr} there would   
an expansion by an isolated $(-3)$--vector, going against Lemma~\ref{l:aux}. 
If $I(S_3)=-3$ then, since $c(S_n)\geq 3$ and by Lemma~\ref{l:n=3} $c(S_3)=1$, 
in Sequence~\eqref{e:contr} there would be two expansions by isolated $(-3)$--vectors, 
one of which would contradict Lemma~\ref{l:extracted}. 
\end{proof}

\begin{prop}\label{p:cleq2(2)}
Let $S\subseteq\bD^n$ be a good subset such that $I(S)+b(S)<0$. Then,
$c(S)\leq 2$.
\end{prop}

\begin{proof}
By definition of bad component, there is a ``minimal'' good subset $S'$ 
associated to $S$ with $I(S')=I(S)$, $b(S')=b(S)$ and $c(S')=c(S)$. 
The set $S'$ is obtained from $S$ via 
contractions by final $(-2)$--vectors, and each bad component 
of $S'$ is of the form $\{v_{s-1},v_s,v_{s+1}\}$, with 
$v_{s-1}\cdot v_{s-1}=v_{s+1}\cdot v_{s+1}=-2$, 
$v_s\cdot v_s < -2$ and $E^{S'}_j=\{s-1,s,s+1\}$ for some $j$ and $s$. 
Then, it is easy to check that if for each bad component of $S'$ we 
apply the transformation
\[
S' \mapsto S'\setminus\{v_s,v_{s+1}\}\cup\{\pi_{e_j}(v_s)\},
\]
the resulting subset $S''$ is good, has no bad components and 
satisfies $c(S'')=c(S)$ and $I(S'')=I(S')+b(S') = I(S) + b(S) <0$. 
Then, by Lemma~\ref{l:cleq2(1)} we have $c(S'')=c(S)\leq 2$. 
\end{proof}

\section{Proof of Theorem~\ref{t:main}: second half}
\label{s:main}

In this section we use the results of Section~\ref{s:prelim} to show that if an oriented 3--manifold $Y$ homeomorphic to a connected sum of lens spaces $L_1\#\cdots\# L_h$ smoothly 
bounds a rational homology 4--ball, then $Y$ is orientation--preserving homeomorphic to 
a connected sum where each summand is (up to orientation) 
homeomorphic to one of the manifolds listed in Theorem~\ref{t:main}(1)--(5). 
The reader is referred to Section~\ref{s:prelim} for the algebraic terminology. 

To start with, we observe that if $h=1$ then it follows from~\cite[Theorem~1.2]{Li}
that $Y=L(p,q)$ with $p/q\in\RR$. Therefore in this case the statement of Theorem~\ref{t:main} is 
established, and from now on we assume $h>1$. 

\begin{figure}[ht]
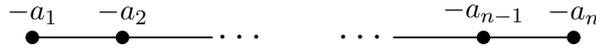

\begin{center}
\setlength{\unitlength}{1mm}
\unitlength=1.2cm
\begin{graph}(8,0.5)(0,0)
\graphnodesize{0.15}
  \roundnode{n1}(1,0)
  \roundnode{n2}(2,0)
  \rectnode{n3}[0,0](3,0)
  \rectnode{n4}[0,0](5,0)
  \roundnode{n5}(6,0)
  \roundnode{n6}(7,0)

  \edge{n1}{n2}
  \edge{n2}{n3}
  \edge{n4}{n5}
  \edge{n5}{n6}

  \autonodetext{n1}[n]{$-a_1$}
  \autonodetext{n2}[n]{$-a_2$}
%  \autonodetext{n3}[n]{$-a_3$}
  \autonodetext{n3}[e]{{\Large$\cdots$}}
  \autonodetext{n4}[w]{{\Large$\cdots$}}
%  \autonodetext{n4}[n]{$-a_{n-2}$}
  \autonodetext{n5}[n]{$-a_{n-1}$}
  \autonodetext{n6}[n]{$-a_n$}
\end{graph}
\end{center}
\caption{The weighted graph prescribing the plumbing $P(p,q)$}
\label{f:fig2}
\end{figure}

Recall that the lens space $L(p,q)$ is orientation--preserving diffeomorphic to the oriented boundary of the 4--dimensional plumbing $P(p,q)$ given by the weighted graph of Figure~\ref{f:fig2}, where $p/q =[a_1,\ldots,a_n]^-$.
It is easy to check that the intersection form of the 4--dimensional plumbing $P(p,q)$ is negative definite. To simplify the notation, we shall also denote by $P_L$ the plumbing associated as above with a lens space $L$. Likewise, if $Y$ is a connected sum of lens spaces $L_1\#\cdots\# L_h$, we define $P_Y$ as the boundary connected sum 
\[
P_{L_1}\natural\cdots\natural P_{L_h}.
\] 
(Note that this is a good definition, i.e. $P_Y$ depends only on $Y$ up to orientation--preserving 
diffeomorphisms). Then, $-Y=(-L_1)\#\cdots\#(-L_h)$. Suppose that $Y=L_1\#\cdots\# L_h$ 
smoothly bounds a rational homology 4--ball $W_Y$. Consider the smooth, closed, negative 
4--manifolds
\[
X_Y := P_Y\cup_{\del} (-W_Y),\quad 
X'_Y := P_{-Y}\cup_{\del} W_Y.
\]
By Donaldson's theorem on the intersection form of definite
4--manifolds~\cite{Do}, the intersection forms of $X_Y$ and $X'_Y$ are
both standard diagonal. 

Suppose that the intersection lattice of $X_Y$ is isomorphic to $\bD^n$ (as defined in Section~\ref{s:prelim}) and the 
intersection lattice of $X'_Y$ is isomorphic to $\bD^{n'}$. Clearly, the groups $H_2(P_Y;\Z)\cong\Z^n$ and 
$H_2(P_{-Y};\Z)\cong\Z^{n'}$ have integral bases which satisfy Equations~\eqref{e:conds}. Therefore, 
via the embeddings $P_Y\subset X_Y$ and $P_{-Y}\subset X'_Y$ we can view such bases as subsets 
$S=S_1\cup\cdots\cup S_h\subset\bD^n$ and $S'=S'_1\cup\cdots\cup S'_h\subset\bD^{n'}$, 
where $S_i$, respectively $S'_i$, originates from a basis  of the corresponding summand $P_{L_i}$ of 
$P_Y$, respectively $P_{-L_i}$ of $P_{-Y}$. Observe that $c(S)=c(S')=h$ and $c(S_i)=c(S'_i)=1$ 
for $i=1,\ldots, h$.

\begin{defn}\label{d:string}
Let $S=\{v_1,\ldots,v_k\}\subset\bD^n$. The {\bf string} of $S$ is $(a_1,\ldots,a_k)$, where 
$a_i=-v_i\cdot v_i$ for $i=1,\ldots, k$. 
\end{defn}

We now briefly recall Riemenschneider's point rule~\cite{Ri}. Let $p>q>0$ be coprime integers, and 
suppose 
\[
\frac pq = [a_1,\ldots,a_l]^-,\ a_i\geq 2,\quad 
\dfrac p{p-q} = [b_1,\ldots,b_m]^-,\ b_j\geq 2.
\]
Then, the coefficients $a_1,\ldots,a_l$ and $b_1,\ldots,b_m$ are related by 
a diagram of the form
\begin{align*}
\bullet \cdots\cdots \bullet\hspace{6cm}\\
\hspace{1.3cm}\bullet \cdots\cdots\bullet \hspace{4.7cm}\\
\hspace{1.8cm}\ddots\hspace{4.2cm} \\
\hspace{3.4cm}\bullet\cdots\cdots \bullet\hspace{2.6cm} \\
\hspace{4.7cm}\bullet \cdots\cdots \bullet\hspace{1.3cm}
\end{align*}
where the $i$--th row contains $a_i-1$ ``points" for $i=1,\ldots,l$, and the first point of each row is 
vertically aligned with the last point of the previous row. The point rule says that there are
$m$ columns, and the $j$--th column contains $b_j-1$ points for every $j=1,\ldots, m$.
For example if $7/5=[2,2,3]^-$ and $7/2 =[4,2]^-$ the corresponding 
diagram is given by 
\begin{align*}
\bullet & \\
\bullet & \\
\bullet & \quad \bullet
\end{align*}

\begin{lem}\label{l:badineq}
We have $b(S)+b(S')\leq c(S)$.
\end{lem}

\begin{proof} 
Suppose that $S_i$ is a bad component of $S$. Then, by Definition~\ref{d:bad} the 
string of $S_i$ is obtained from the string $(2,a,2)$, $a\geq 3$, via a finite sequence of operations of 
type Lemma~\ref{l:fraction}(1) or Lemma~\ref{l:fraction}(2). Applying Riemenschneider's point 
rule (cf.~\cite[Proof of Lemma~2.6]{Li}) one can easily check that, similarly, the string of  $S'_i$ is 
obtained from the string
\[
(3,\overbrace{2,\ldots,2}^{a-3},3)
\]
via a sequence of operations of type Lemma~\ref{l:fraction}(1) or Lemma~\ref{l:fraction}(2). 
Clearly, such a string is not the string of a bad component. This shows that, for each $i=1,\ldots, h$, at 
most 
one set among $S_i$ and $S'_i$ can be a bad component, which implies the statement of the lemma. 
\end{proof}

\begin{lem}\label{l:ineq}
Up to replacing $S$ with $S'$ we have $I(S)+c(S)\leq 0$ and
$I(S)+b(S)<0$. 
\end{lem}

\begin{proof}
We have $c(S)=c(S')$ and by Lemma~\ref{l:negsum} it follows that
$I(S)+I(S')=-2c(S)$. Therefore, up to replacing $S$ with $S'$ we may assume 
$I(S)\leq -c(S)$. If $I(S)<-c(S)$ then, since $b(S)\leq c(S)$, we have 
\begin{equation}\label{e:1}
I(S) + b(S) \leq I(S) + c(S) < 0
\end{equation}
and the lemma holds. Now assume $I(S)=-c(S)=I(S')$. By Lemma~\ref{l:badineq}
$b(S)+b(S')\leq c(S)$. Therefore, up to replacing $S$ with $S'$ we may assume $I(S)=-c(S)$ 
and $b(S)<c(S)$, so it follows that 
\begin{equation}\label{e:2} 
I(S) + b(S) < I(S) + c(S) = 0.
\end{equation}
This concludes the proof.
\end{proof}

In view of Lemma~\ref{l:ineq}, from now on we will assume: 
\[
I(S)+c(S)\leq 0\quad\text{and}\quad I(S)+b(S)<0. 
\]

Now there are two possibilities. Either $S$ is irreducible or it is reducible. 
We consider the two cases separately. 

\part*{First case: $S$ irreducible} 

In this case $S$ is a good subset by definition. Since $I(S)+b(S)<0$, by 
Proposition~\ref{p:cleq2(2)}, we have $h=c(S)\leq 2$. But we are assuming $h>1$, so we 
conclude $h=c(S)=2$. Since $b(S)\leq c(S)$, we must analyze the three subcases $b(S)=0$, 
$b(S)=1$ and $b(S)=2$. 

{\em First subcase: $b(S)=0$}. By Proposition~\ref{p:I<0} there is a sequence of 
contractions of good subsets without bad components as in~\eqref{e:expans}. 
Since $I(S_3)\leq I(S)\leq -c(S)=-2$, applying~\eqref{e:Ic1}, \eqref{e:Ic2} and 
Lemma~\ref{l:n=3} it is easy to check that  $(I(S),c(S))=(-2,2)$ and 
either $(I(S_3),c(S_3))=(-2,2)$ or $(I(S_3),c(S_3))=(-3,1)$. The latter case 
is excluded by~\eqref{e:Ic1},~\eqref{e:Ic2} and  Lemma~\ref{l:extracted}. 
So we are left with the case $(I(S_3),c(S_3))=(-2,2)$. 
Since $I(S)=I(S_3)=-2$, by~\eqref{e:Ic1} and~\eqref{e:Ic2} 
Sequence~\eqref{e:expans} satisfies the assumptions of Lemma~\ref{l:I=-2}. 
Applying the lemma to $S_k$ for each $k=3,\ldots,n$ it is easy to check that 
the string associated with $S$ is a union of two  strings $s_1\cup s_2$ related 
to each other by Riemenschneider's point rule. We conclude that, up 
to orientation, $Y=L_1\# L_2$ with $L_1=L(p,q)$ and $L_2=L(p,p-q)$ for some 
$p>q>0$, so Case (2) of Theorem~\ref{t:main} holds.

{\em Second subcase: $b(S)=1$}. By definition of bad component there is a 
sequence of contractions of good subsets 
\[ 
S\searrow\cdots\searrow T
\] 
obtained by erasing final $(-2)$--vectors, such that $c(T)=c(S)=2$, $b(T)=b(S)=1$, 
$I(T)=I(S)\leq -c(S)=-2$ and such that the bad component of $T$ is of the 
form $\{v_{s-1}, v_s, v_{s+1}\}$ with $v_{s-1}\cdot v_{s-1}=v_{s+1}\cdot v_{s+1}=-2$ 
and $v_s\cdot v_s<-2$. We can further reduce $T$ to a good subset without 
bad components $U$ as in the proof of Proposition~\ref{p:cleq2(2)}, so that 
$I(U)=I(T)+1\leq -1$, $c(U)=c(T)=2$ and one connected component of the 
intersection graph of $U$ consists of a single element. Moreover, by 
Proposition~\ref{p:I<0} there is a sequence of contractions 
\[
U\searrow\cdots\searrow S_3\subset\bD^3
\] 
satisfying~\eqref{e:Ic1} and~\eqref{e:Ic2}. In particular, 
$I(U)\geq I(S_3)\geq -3$. The case $I(U)=I(S_3)=-3$ is impossible 
by~\eqref{e:Ic1} and~\eqref{e:Ic2} because $c(U)=2$ while by Lemma~\ref{l:n=3} 
the equality $I(S_3)=-3$ implies $c(S_3)=1$.

If $I(U)=-2$, the same analysis made in the first subcase above
shows that $(I(S_3),c(S_3))=(-2,2)$ and the string of $U$ must be of the form 
\[
(n)\cup (\overbrace{2,\cdots, 2}^{n-1}),\ n\geq 2,
\]
therefore the string of $T$ is of the form 
\[
(2,n+1,2)\cup (\overbrace{2,\cdots, 2}^{n-1}),\quad n\geq 2.
\]
We conclude that the string of $S$ is of the form 
\[
s_n \cup (\overbrace{2,\cdots, 2}^{n-1}),\quad n\geq 2, 
\]
where the string $s_n$ is obtained from $(2,n+1,2)$ by a finite sequence of 
operations as in Lemma~\ref{l:fraction}. Applying Lemma~\ref{l:fraction} 
we conclude that, up to orientation, $Y=L_1\# L_2$ with $L_1=L(m^2n, mnk+1)$, 
$L_2=L(n,n-1)$ for some $n\geq 2$, $m>k>0$ and $(m,k)=1$. In other words,  
Case (4) of Theorem~\ref{t:main} holds.

If $I(U)=-1$ then $I(S)=I(T)=I(U)-1=-2$. By Lemma~\ref{l:negsum} the subset 
$S'\subset\bD^{n'}$ originating from the integral basis of $H_2(P_{-Y};\Z)$ 
satisfies $(I(S'),c(S'))=(-2,2)$ and by Lemma~\ref{l:badineq} we have $b(S')\leq 1$.
If $b(S')=0$, the argument given in the first subscase shows that 
$Y=L(p,q)\# L(p,p-q)$, so Case (2) of Theorem~\ref{t:main} holds. If $b(S')=1$, the same 
argument just used for $S$ shows that there is a 
sequence of contractions of good subsets 
\[ 
S'\searrow\cdots\searrow T'
\] 
obtained by erasing final $(-2)$--vectors, such that $c(T')=c(S')=2$, 
$b(T')=b(S')=1$, $I(T')=I(S')=-2$. By Riemenschneider's point rule, 
the string of the bad component of $T'$ is dual to the string of the non--bad 
component of $T$, and vice--versa. This immediately implies that the intersection graph of $S$ is a 
union of connected components $S_1\cup S_2$, where (i) the string of $S_1$ is obtained from 
$(2,a,2)$, for some $a\geq 3$, by a finite sequence of operations as in Lemma~\ref{l:fraction} and 
(ii) the string of $S_2$ is the Riemenschneider dual of a string obtained in a similar way from 
$(2,b,2)$, for some $b\geq 3$. Clearly $I(S_1)=a-5$ and by Lemma~\ref{l:negsum} we have 
$I(S_2)=-b+3$. But $I(S)=I(S_1)+I(S_2)=-2$ implies $a=b$, therefore applying 
Lemma~\ref{l:fraction} we conclude that, up to orientation, $Y=L_1\# L_2$ with $L_1, 
-L_2\in\FF_{a-1}$. 
Hence, Case (5) of Theorem~\ref{t:main} holds.

{\em Third subcase: $b(S)=2$}. 

As in the previous subcase, there is a sequence of contractions of good subsets by 
final $(-2)$--vectors 
\[
S\searrow\cdots\searrow T
\]
such that $c(T)=c(S)=2$, $b(T)=2$ and $I(T)=I(S)$. Erasing final $(-2)$--vectors as before  
we get a good subset $U\subseteq\bD^2$, and it is easy to check that the string of $U$ 
must be equal to $(2)\cup (2)$. Therefore the string of $T$ is  
\[
(2,3,2)\cup (2,3,2)
\]
and the string of $S$ is of the form $s_n\cup t_n$, where each one of $s_n$ and $t_n$ is obtained 
from $(2,3,2)$ by a finite sequence of operations as in Lemma~\ref{l:fraction}. 
Applying Lemma~\ref{l:fraction} as in the previous subcase we see that, up to 
orientation, $Y=L_1\# L_2$ with $L_1=L(2m^2,2mk+1)$ and $L_2=L(2p^2,2ph+1)$, for some 
$m>k>0$, $(m,k)=1$ and $p>h>0$, $(p,h)=1$. Therefore Case (3) of Theorem~\ref{t:main} 
holds. Summarizing, so far we have proved:
\begin{lem}\label{l:sofar}
If $h\geq 2$ and the subset $S$ is irreducible, then $h=2$ and $Y$ is (possibly 
after reversing its orientation) homeomorphic to one of the manifolds given in 
Theorem~\ref{t:main}(2)--(5).
\qed\end{lem}

\part*{Second case: $S$ reducible} 

In this case the set $S$ can be written as a disjoint union $S=\cup_j T_j$ 
of maximal irreducible subsets $T_j\subset S$. As observed in~\cite[Remark~2.1]{Li}, 
the elements of $S$ are linearly independent over $\Z$ because they satisfy 
Equation~\eqref{e:conds}. We claim that, for each index $j$, 
\begin{equation}\label{e:irredcond}
|T_j| = |\cup_{v_i\in T_j} V_i|.
\end{equation}
In fact, by the linear independence we have 
\[
|\cup_{v_i\in T_j} V_i|\geq |T_j|,
\]
while by the maximality of each $T_j$ the union $\bigcup_j\left(\cup_{v_i\in T_j} V_i\right)$ 
is disjoint. Hence 
\[
|\bigcup_j\left(\cup_{v_i\in T_j} V_i\right)| = \sum_j |\cup_{v_i\in T_j} V_i|.
\]
So we have 
\[
\sum_j |\cup_{v_i\in T_j} V_i| \geq \sum_j |T_j| = |S| 
\geq |\bigcup_j\left(\cup_{v_i\in T_j} V_i\right)| = \sum_j |\cup_{v_i\in T_j} V_i|, 
\]
and the claim follows. Next, we claim that there exists at least one index $i$ 
such that $I(T_i)+c(T_i)\leq 0$ and $I(T_i)+b(T_i)<0$. Since $b(T_j)\leq c(T_j)$
for every $j$, it suffices to show that either (i) there is an index $i$ such that 
$I(T_i) + c(T_i)<0$ or (ii) there is an index $i$ such that $I(T_i)+c(T_i)=0$ and 
$I(T_i)+b(T_i)<0$. Since 
\begin{equation*}
\sum_j (I(T_j) + c(T_j)) = I(S)+c(S)\leq 0, 
\end{equation*}
if (i) does not hold then $I(T_j)+c(T_j)=0$ for 
every $j$. In this case, since we also have 
\[
\sum_j (I(T_j) + b(T_j)) = I(S)+b(S) < 0, 
\]
we conclude that $I(T_i) + b(T_i)<0$ for some $i$, i.e.~that (ii) holds. By~\eqref{e:irredcond} and 
the last claim, we may apply Lemma~\ref{l:sofar} to $T_i$, viewed as a subset of the span of the 
$e_k$'s hit by the vectors of $T_i$. Therefore, after possibly renaming the $L_i$'s we conclude that 
either $L_h$ or $L_{h-1}\# L_h$ smoothly bounds a rational homology ball $Z$. Let $p\in\del Z$, 
let $B\subseteq \del Z$ a regular neighborhood of $p$, and let $W$ be a smooth rational homology ball 
with boundary $L_1\#\cdots\# L_h$. Clearly, the 3--manifold with boundary $\del Z\setminus B$ can be viewed 
as a subset of $\del W$ as well as $\del Z$. Then, the space  
\[
W\cup_{(\del Z)\setminus B} (-Z)
\] 
obtained by gluing $W$ and $-Z$ together along $\del Z\setminus B$, 
is (after smoothing corners)  a smooth rational homology ball with boundary $L_1\#\cdots\# L_{h'}$, where
$h'$ is equal to $h-1$ or $h-2$, respectively. Thus, we have proved:
\begin{lem}\label{l:final} 
If $h\geq 2$ and $S$ is reducible, after possibly renaming the $L_i$'s one of the following 
holds:
\begin{itemize}
\item
$L_1\#\cdots\# L_{h-1}$ and $L_h$ smoothly bound a rational homology ball.
\item
$L_1\#\cdots\# L_{h-2}$ and $L_{h-1}\# L_h$ smoothly bound a rational homology ball.
\qed
\end{itemize}
\end{lem}
At the beginning of this section we observed that when $h=1$ the statement of Theorem~\ref{t:main} was proved in~\cite{Li}. Therefore, the second half of Theorem~\ref{t:main} follows combining   
Lemmas~\ref{l:sofar} and~\ref{l:final}.

\end{document}